\documentclass{amsart}

\usepackage{amssymb}

\newtheorem{theorem}{Theorem}[section]
\newtheorem{proposition}[theorem]{Proposition}
\newtheorem{conjecture}[theorem]{Conjecture}
\newtheorem{lemma}[theorem]{Lemma}
\newtheorem{remark}[theorem]{Remark}

\renewcommand{\P}{{\mathbb P}}
\newcommand{\C}{{\mathbb C}}
\newcommand{\Z}{{\mathbb Z}}
\newcommand{\Q}{{\mathbb Q}}
\newcommand{\R}{{\mathbb R}}
\newcommand{\X}{{\mathbb X}} 

\newcommand{\by}{{\mathcal Y}} 
\newcommand{\bl}{{\mathcal L}} 
\newcommand{\bz}{{\mathcal Z}} 
\newcommand{\bigx}{{\mathcal X}}

\newcommand{\aut}{\mathop{\rm Aut}\nolimits}
\newcommand{\pic}{\mathop{\rm Pic}\nolimits}
\renewcommand{\phi}{\varphi}

\begin{document}
 
\title{On an example of Aspinwall and Morrison} 

\author{Bal\' azs Szendr\H oi}

\thanks{Research partially supported by an Eastern European Research Bursary from Trinity College, Cambridge and an ORS Award from the British Government.}

\address{Mathematics Institute, University of Warwick, Coventry CV4 7AL, United Kingdom}

\email{balazs@maths.warwick.ac.uk}

\date{23 June 2002}


\begin{abstract} 
In this paper, a family of smooth multiply connected Calabi--Yau threefolds 
is investigated. The family presents a counterexample to global Torelli
as conjectured by Aspinwall and Morrison.
\end{abstract}

\maketitle
 
\section*{Introduction}

The aim of this paper is to prove 

\begin{theorem} The one-parameter family of smooth, multiply connected 
Calabi--Yau threefolds $\by\rightarrow B$ over the base 
$B=\P^1 \setminus\{1,\xi, \ldots, \xi^4,\infty\}$, 
constructed by Aspinwall--Morrison in~\cite{aspinwall_morrison} 
(cf. Section~\ref{am_defs}), with $\xi$ a primitive fifth root of unity, 
has the following properties:
\begin{itemize} 
\item For any $t\in B$, there exists an isomorphism 
\[ H^3 (Y_t, \Q) \cong H^3 (Y_{\xi t}, \Q)\] 
preserving rational polarized Hodge structures (for a stronger statement, 
see Theorem~\ref{hodgsame}).
\item There is a Zariski-open set $U\subset B$ such that for $t\in U$, 
$i=0,\ldots ,4$, the fibres $Y_{\xi^it}$ are pairwise non-isomorphic as 
algebraic varieties.
\end{itemize}\label{am_main_theorem}\end{theorem}

The family $\by\rightarrow B$ is a quotient of a family of quintics, 
manufactured in such a way that a certain symmetry of a cover 
$\bz\rightarrow B$ of $\by\rightarrow B$ fails to descend in any obvious way 
to a symmetry of $\by\rightarrow B$. The existence of this symmetry on the 
cover implies the statement about Hodge structures (Theorem~\ref{hodgsame}). 
On the other hand, an isomorphism between $Y_t$ and $Y_{\xi t}$ for general 
$t$ would force, via a specialization argument 
(Theorem~\ref{am_onestep_theorem}), the existence of an automorphism $\sigma$ 
on the fibre $Y_0$ over $0$ of a special kind. However, the automorphism group
of $Y_0$ can be computed explicitly (Theorem~\ref{autY}), and such a $\sigma$ 
does not exist. For technical reasons, the argument runs on a family of
singular models $\bar\by\rightarrow B$ of $\by\rightarrow B$. 
(See Section 4.) 

Theorem~\ref{am_main_theorem} establishes the fact, conjectured by Aspinwall 
and Morrison, that the family $\by\rightarrow B$ provides a counterexample to 
global Torelli for Calabi--Yau threefolds. Previous counterexamples 
to Torelli were given in~\cite{tcex}; there I considered families of
birationally equivalent Calabi--Yau threefolds. By \cite[Theorem 4.12]{flops}, 
birational equivalence implies isomorphism between (rational) Hodge 
structures. However, in the present case the situation should be entirely 
different. 

\begin{conjecture} For general $t\in B$, the threefolds $Y_{\xi^it}$ for 
$i=0,\ldots, 4$ are not birationally equivalent to one another. 
\label{notbir}
\end{conjecture} 

One obvious direct approach to this conjecture is to aim to understand the 
various birational models of a fixed fibre $Y_t$. Birational models
of minimal threefolds can be studied via their cones of nef divisors in the 
Picard group; so this method requires an explicit understanding of
the nef cone of $Y_t$. An \'etale cover $Z_t$ of $Y_t$ is a toric 
hypersurface. A recent conjecture \cite[Conjecture 6.2.8]{cox_katz} of Cox 
and Katz aimed at giving a complete understanding of the nef cone of toric
Calabi--Yau hypersurfaces. However, I prove in~\cite{counterex} that in fact 
the conjecture of Cox and Katz fails for $Z_t$. At this point the computation 
of the nef cone of $Y_t$ seems rather hopeless. A different approach to 
Conjecture~\ref{notbir} is required.

To conclude the introduction, let me point out that the varieties $Y_t$ are 
multiply connected with fundamental group $\Z/5\Z$ (Proposition~\ref{extends} 
and Proposition~\ref{fundZ}). This is a curious fact. The construction of 
Aspinwall and Morrison requires in an essential way that 
members of the mirror Calabi--Yau family should have nontrivial (and in fact 
non-cyclic) fundamental group. Computations of Gross~\cite[Section 3]{gross} 
connect torsion in the integral cohomologies of mirror Calabi--Yau threefolds,
and these computations imply that the cohomology (and hence homology) of $Y_t$
should have torsion of some kind. However, the direct relationship between 
failure of Torelli and the fundamental group seems rather mysterious; compare 
also Remark~\ref{nonact}. 

\vspace{0.1in}
 
\noindent
{\bf Acknowledgments} \ I wish to thank Pelham Wilson, Mark Gross and Peter 
Newstead for comments and help. 

\vspace{0.1in}
 
\noindent
{\bf Notation and conventions} \ All schemes and varieties are defined over 
$\C$. A {\it Calabi--Yau threefold} is a normal projective threefold $X$ with 
canonical Gorenstein singularities, satisfying $K_X\sim 0$ and 
$H^1(X,{\mathcal O}_X)=0$. Some statements use the language of 
toric geometry; my notation follows Fulton~\cite{fulton} and 
Cox--Katz~\cite[Chapter~3]{cox_katz}. If $A$ is a $\Z$-module then 
$A_{\rm free}$ denotes the torsion free part. 

\section{The construction} 
\label{am_defs}

Following~\cite{aspinwall_morrison}, define maps $g_i:\P^4 \rightarrow \P^4$ by
\setlength{\arraycolsep}{2pt}
\[
\begin{array}{rccl}
g_1: & [z_0 : z_1 : z_2 : z_3 : z_4] & \mapsto &  [z_0 : \xi z_1 : \xi^2 z_2 : \xi^3 z_3 : \xi ^4 z_4]\\
g_2: & [z_0 : z_1 : z_2 : z_3 : z_4] & \mapsto &  [z_0 : \xi z_1 : \xi ^3 z_2 : \xi z_3 : z_4] \\
g_3: & [z_0 : z_1 : z_2 : z_3 : z_4] & \mapsto &  [z_1 : z_2 : z_3 : z_4 : z_0]

\end{array}
\]
where $\xi$ is a fixed primitive fifth root of unity. Let
\[G=\langle g_1, g_2, g_3 \rangle,\,\, H=\langle g_1, g_2 \rangle\]
be subgroups of $PGL(5,\C)$ generated by the transformations $g_i$.
As abstract groups $H\cong (\Z/5\Z)^2$, $G \cong \Z/5\Z \ltimes (\Z/5\Z)^2$.

I will be interested in hypersurfaces in the varieties $\P^4/G$ and 
$\P^4/H$; the latter is a toric variety and its toric description will be 
useful in the sequel. 

\begin{proposition} In the contravariant description, 
$\P^4/H\cong\P_{M,\Delta}$, where 
$M\cong\Z^4$ and $\Delta\subset M_\R$ is the polyhedron
\[ \Delta = {\mathop{\rm span}}\{ (1,0,0,0), (-3,5,-4,-2), (0,0,1,0), (0,0,0,1), (2,-5,3,1)\}.\]
With $N={\rm Hom}(M,\Z)$, the dual polyhedron $\Delta^*\subset N_\R$ of 
$\Delta$ is 
\[\Delta^*={\mathop{\rm span}}\{(-1,-2,-1,-1), (4,1,-1,-1), (-1,-1,-1,-1),
(-1,2,4,-1), (-1,0,-1,4)\}.\]

\noindent
The polyhedron $\Delta^*$ has no interior lattice points apart from the
origin, has no lattice points in the interiors of its three- or
one-dimensional faces, and has precisely two lattice points 
$P_{2i-1}, P_{2i}$, $i=1,\ldots, 10$ in the interiors of each of its ten 
two-dimensional faces. 
\end{proposition} 
\begin{proof}
This is a standard toric calculation; for details see~\cite[Proposition 1.1]{counterex}. 
\end{proof} 

Let $\Sigma$ be the fan consisting of cones over faces of $\Delta^*$ 
in $N_\R$. This fan defines the toric variety 
$\X_{N,\Sigma}\cong \P_{M,\Delta}$. 

\begin{proposition} $\P_{M,\Delta}$ is a $\Q$-factorial Gorenstein variety, 
with ten curves of ca\-nonical singularities. Every permutation 
$\eta$ of the lattice points $\left\{ P_i\right\}$ gives rise to 
a partial resolution $\X_{\Sigma_\eta}\rightarrow\P_{M,\Delta}$.
The varieties $\X_{\Sigma_\eta}$ have isolated singularities only.
\end{proposition}
\begin{proof}
This is basic toric geometry. The curves of singularities correspond to 
the ten two-dimensional faces of $\Delta^*$. The singularities can be 
partially resolved by subdividing the fan $\Sigma$ using the lattice points 
$\left\{ P_i\right\}$ in any order. Any permutation $\eta$ of these points  
gives a fan $\Sigma_\eta$ in the space $N_\R$ and a corresponding 
toric partial resolution $\X_{\Sigma_\eta}$ with isolated singularities.  
\end{proof}

The family of hypersurfaces of interest in this paper is constructed from 
\[ {\mathcal Q} = \left\{ \sum_{i=0}^4 z_i^5 -5t \prod_{i=0}^4 z_i =0 \right\}\subset \P^4\times B, \]
where $B=\C\setminus \left\{1,\xi, \ldots, \xi^4\right\}$. Second projection 
gives a smooth family $p:{\mathcal Q}\rightarrow B$ of 
Calabi--Yau quintics~$Q_t$. The groups $G$ and $H$ act on $\P^4\times B$ 
by acting trivially on $B$, and hence on ${\mathcal Q}$; 
these actions preserve holomorphic three-forms in the fibres. Let
\[
\begin{array}{rcccl}
\bar {\mathcal Z} & = & {\mathcal Q} / H, \\
\bar {\mathcal Y} & = & {\mathcal Q} / G & = &  \bar {\mathcal Z}/K.
\end{array}
\]
Here $K\cong \Z/5\Z$ is the group generated by the image of $g_3$ in 
$\aut(\bar {\mathcal Z})$. Both $\bar \bz$ and $\bar \by$ are naturally
families over $B$ with fibres $\bar Z_t$ and $\bar Y_t$ respectively. 

\begin{proposition} For $t\in B$, $\bar Z_t$ is a canonical Calabi--Yau 
threefold with ten isolated $\frac{1}{5}(1,1,3)$ quotient singularities. 
The group $K$ acts freely on $\bar Z_t$. The variety $\bar Y_t$ is a 
canonical Calabi--Yau threefold\ with two isolated $\frac{1}{5}(1,1,3)$ 
quotient singularities.
\end{proposition} 
\begin{proof} Easy explicit check.    
\end{proof}

The family $\bar \bz\rightarrow B$ is a family of non-degenerate anti-canonical 
hypersurfaces in the toric variety $\P_\Delta$. The partial resolutions
$\X_{\Sigma_\eta}\rightarrow\P_{M,\Delta}$ give rise to morphisms
$\bz_{\eta}\rightarrow \bar \bz$ over $B$, with $\bz_{\eta}\rightarrow B$ 
a family of nonsingular threefolds as $\X_{\Sigma_\eta}$ is nonsingular 
in codimension three.

\begin{proposition} The families $\bz_{\eta}$ are all canonically isomorphic 
to a unique toric resolution $\bz\rightarrow \bar \bz$ over $B$. For $t\in B$,
the fibre $Z_t$ is a smooth Calabi--Yau threefold with Hodge numbers 
$h^{1,1}(Z_t)=21$, $h^{2,1}(Z_t)=1$. In the resolution 
$Z_t\rightarrow \bar Z_t$ there are two exceptional divisors over every 
singular point $S_i$, a Hirzebruch surface $E_i\cong{\mathbb F}_3$ and a 
projective plane $F_i\cong\P^2$ intersecting in a $\P^1$ which is the negative
section in the Hirzebruch surface and a line in $\P^2$. 
\label{unique}
\end{proposition}
\begin{proof} Let $\eta_1$, $\eta_2$ be two permutations of the interior 
lattice points. There is a corresponding birational map 
$\X_{\Sigma_{\eta_1}} \dashrightarrow \X_{\Sigma_{\eta_2}}$ whose exceptional 
sets are disjoint from the families $\bz_{\eta_i}$. This implies the first 
part. The other statements follow from easy toric calculations. 
\end{proof}

\begin{proposition} The action of the group $K\cong \Z/5\Z$ on $\bar\bz$ 
extends to a free action on the resolution $\bz$ over $B$. Thus there is an 
\'etale cover $\bz \rightarrow \by=\bz/K$ over $B$. The fibre $Y_t$ for 
$t\in B$ is a Calabi--Yau resolution of $\bar Y_t$ with Hodge numbers 
$h^{1,1}(Y_t)=5$, $h^{2,1}(Y_t)=1$. 
\label{extends}
\end{proposition}
\begin{proof} The action of $K$ is generated by the symmetry $g_3$ of $\P^4$. This 
symmetry descends to the toric variety $\P_\Delta$ as a toric symmetry 
induced by a lattice isomorphism $\alpha_3: M\rightarrow M$ fixing the 
polyhedron $\Delta$ and permuting the lattice points $\{P_i\}$. 
Composition with the permutation induced by $\alpha_3$ gives a correspondence 
$\eta\rightarrow \eta'$ between permutations of the set $\{P_i\}$ and 
$\alpha_3$ gives rise to an isomorphism 
$\tilde g_3: \X_{\Sigma_\eta}\rightarrow \X_{\Sigma_{\eta'}}$. 
This isomorphism restricts to anti-canonical families as an isomorphism 
$\bz_{\eta}\rightarrow \bz_{\eta'}$, or, by Proposition~\ref{unique}, 
as an automorphism $\bz\rightarrow \bz$. By construction, this automorphism
is the required extension of $g_3$ and it clearly generates a free group 
action on $\bz$ over $B$.
\end{proof}

I conclude this section by proving two auxiliary statements. 

\begin{proposition} The family $\bar\by\rightarrow B$ restricted to 
a neighbourhood of $0\in B$ is the universal 
deformation space of its central fibre $\bar Y_0$ in the analytic category.
\label{betterprop}
\end{proposition}
\begin{proof} By general theory, the projective variety $\bar Y_0$ 
has a versal deformation space $\bigx\rightarrow S$ 
in the analytic category. $\bar Y_0$ is a canonical Calabi--Yau threefold.  
Thus $H^0(\bar Y_0, T_{\bar Y_0})=0$ and this implies 
that $\bigx\rightarrow S$ is in fact universal. 
By Ran's extension~\cite{ran_klein} of the Bogomolov--\-Tian--\-Todorov 
theorem, unobstructedness holds for $\bar Y_0$. Thus $S$ is smooth.
Further, the codimension of the singularities of $\bar Y_0$ is three.  
By the argument of~\cite[A.4.2]{cox_katz}, it follows that
the first-order tangent space of $S$ at the base point 
is isomorphic to $H^1(\bar Y_0, T_{\bar Y_0})$, a one-dimensional 
complex vector space. 

In order to prove that $\bar\by\rightarrow B$ is the universal deformation 
space, all I need to show is that its Kodaira-Spencer map is injective. 
Recall the family ${\mathcal Q}\rightarrow B$, a deformation of the Fermat 
quintic $Q_0$ over $B$. Choosing a ($G$-invariant) three-form on $Q_0$ gives
rise to a commutative diagram
\[\begin{array}{ccccc}
T_0(B)&\stackrel{k}{\longrightarrow} &H^1(Q_0, T_{Q_0})&\stackrel{\sim}{\longrightarrow} &H^1(Q_0,\Omega^2_{Q_0})\\
\Vert && && \Big\uparrow\rlap{$\vcenter{\hbox{$\scriptstyle j$}}$}\\
T_0(B)&\stackrel{l}{\longrightarrow} &H^1(\bar Y_0, T_{\bar Y_0})&\stackrel{\sim}{\longrightarrow} &H^1(\bar Y_0,\hat \Omega^2_{\bar Y_0}).
\end{array}\]
Here $k$ and $l$ are the Kodaira--Spencer maps, whereas the map
$j$ is given by pullback of (orbifold) 
two-forms (the sheaf of orbifold two-forms $\hat \Omega^2_{\bar Y_0}$ 
is defined carefully in~\cite[A.3]{cox_katz}). 
The map $k$ is injective, as ${\mathcal Q}$ is a nontrivial 
first-order deformation of $Q_0$. By commutativity, $l$ is also injective. 
This proves the Proposition.
\end{proof}

\begin{proposition} For $t\in B$, the Calabi--Yau manifold $Z_t$ is 
simply connected.
\label{fundZ}
\end{proposition}
\begin{proof} The variety $Z_t$ is a resolution of the threefold
$\bar Z_t= Q_t/H$. Let $Q_t^0$ be the open set of $Q_t$ on which 
the action of $H$ is free; it is the complement of a 
finite set of points and hence is simply connected.
Let $Z_t^0=Q_t^0/H$; $\pi_1(Z_t^0)\cong H$. 

The fundamental group of $Z_t$ is a quotient group of $H$. 
Let $T_t$ be the universal cover of $Z_t$; by the generalized 
Riemann existence theorem, $T_t$ is an algebraic variety and it clearly has
trivial canonical bundle. Let $T_t^0$ be the preimage of $Z_t^0$ under the 
covering map. Then $T_t^0$ has finite fundamental group; 
let $\tilde T_t^0$ be its universal cover. 
$\tilde T_t^0$ is an algebraic variety again. Notice however, 
that $Q_t^0$, $\tilde T_t^0$ are both universal covers of the variety 
$Z_t^0$, and thus by the uniqueness part of the generalized Riemann 
existence theorem they must be isomorphic. Thus there exists a diagram 
\[
\begin{array}{ccccc}
Q_t &\supset& Q_t^0  \\
&&            \downarrow\\
\Big\downarrow && T_t^0 & \subset & T_t \\ 
&&    \downarrow && \downarrow \\
\bar Z_t & \supset & Z_t^0 & \subset & Z_t.  
\end{array}
\]
The covering $Q_t^0\rightarrow T_t^0$ corresponds to a group $L$ of 
holomorphic 
automorphisms of $Q_t^0$. An automorphism of $Q_t^0$ can be thought of as a 
birational self-map of $Q_t$. However, as 
$Q_t$ is a minimal Calabi--Yau threefold 
with Picard number one, it has no birational self-maps with a 
nontrivial exceptional locus. So $L$ consists of 
automorphisms of $Q_t$. The fact that the map $Q_t^0\rightarrow T_t^0$ factors 
the map $Q_t^0\rightarrow Z_t^0$ implies that $L$ must be a subgroup of $H$.

Thus I conclude that $T_t$ is birational to a quotient $Q_t/L$ for 
a subgroup $L$ of $H$. 
Moreover, $\chi(Z_t)=40$ so $\chi(T_t)$ equals either 40, 200 or 1000. 
On the other hand, for every subgroup $L$ of $H$, the quotient $Q_t/L$ has
a Calabi--Yau desingularization. 
As the Euler number is a birational invariant of smooth Calabi--Yau threefolds,
the Euler number of this desingularization must be equal to that of $T_t$. 
Finally, it is easy to check that $H$ has no subgroup 
$L$ such that a Calabi--Yau desingularization of $Q_t/L$ has Euler number 
200 or 1000. Thus $L=H$ and so $T_t=Z_t$ is its own universal cover.  
\end{proof}

\section{Hodge structures} \label{hodge}

Let  $Z$, $Y$ denote the differentiable manifolds underlying
the fibres $Z_t$, $Y_t$. Let $V_Z=H^3(Z, \Z)_{\rm free}$, $V_Y=H^3(Y, \Z)_{\rm free}$,
with antisymmetric pairings $Q_Z$, $Q_Y$ given by cup product. 

\begin{lemma} \label{modules}
Pullback by the map $\pi:Z\rightarrow Y$ induces an injection 
\[ \pi^*: V_Y \hookrightarrow V_Z
\]
with image of index at most 25. Under this embedding, 
\[ Q_Z(\pi^*\alpha_1, \pi^*\alpha_2) = 5 \, Q_Y(\alpha_1, \alpha_2).
\]
Consequently, there is an embedding of groups
\[ \aut_\Z(V_Z, Q_Z)\stackrel{j}{\longrightarrow}\aut_\Q(V_Y\otimes \Q, Q_Y).\]
\end{lemma}
\begin{proof}
The group $K\cong \Z/5\Z$ acts without fixed points on $Z$, so
the map $\pi$ induces a spectral sequence
\[ E_2^{p,q}=H^p(K; H^q(Z, \Z)) \Rightarrow H^{p+q}(Y, \Z). 
\]
The terms $E_2^{p.q}$ for $p>0$ are torsion, so $V_Y=(E_\infty^{0,3})_{\rm free}$.  
On the other hand, $(E_2^{0,3})_{\rm free}=H^0(K, H^3(Z,\Z)_{\rm free})=(V_Z)^K$. 
There are two differentials from $(E_2^{0,3})$,
both having image $\Z/5\Z$.  So there is an
injection
\[ \pi^*:V_Y \hookrightarrow (V_Z)^K
\]
with image of index at most 25. This map is an isomorphism 
when tensored~by~$\Q$. As both $V_Z$ and $V_Y$ have rank 
four, $K$ must act trivially on $V_Z$ and this proves the first part.
The other two statements are immediate.
\end{proof}

Let ${\mathcal D}_Y$ be the period domain parameterizing weight 3 polarized 
Hodge structures on $(V_Y, Q_Y)$. Fixing a point $t\in B$, a marking 
$H^3(Y_t, \Z)_{\rm free}\cong V_Y$ and a universal cover 
$\tilde B$ of $B$ leads to holomorphic period maps
\[\begin{array}{ccc} 
\tilde B & \stackrel{\tilde\psi}{\longrightarrow} & {\mathcal D}_Y \\
\downarrow && \downarrow \\
B & \stackrel{\psi}{\longrightarrow} & {\mathcal D}_Y/\Gamma
\end{array} 
\]
where $\Gamma$ is any subgroup of $\aut_\Q(V_Y\otimes \Q, Q_Y)$ 
containing all geometric monodromies and acting properly discontinuously on 
${\mathcal D}$.  
Choose \[\Gamma=j(\aut_\Z(V_Z, Q_Z))\subset \aut_\Q(V_Y\otimes \Q, Q_Y)\]
under the embedding $j$ of Lemma~\ref{modules}.  

\begin{lemma} $\Gamma$ acts properly discontinuously on ${\mathcal D}_Y$, so 
${\mathcal D}_Y/\Gamma$ is an analytic space.
\end{lemma}
\begin{proof} See~\cite[Section I.2]{grif}.\end{proof} 

After all these preparations, I can state 

\begin{theorem} For $\Gamma$ chosen as above,  
the period map $\psi: B\rightarrow{\mathcal D}_Y/\Gamma$
is of degree at least five. More precisely, 
if $t_1, t_2\in B$ satisfy $t_1^5=t_2^5$, then 
$\psi(t_1) = \psi(t_2)$. In particular, $Y_{t_1}$ and $Y_{t_2}$ have 
isomorphic rational Hodge structure. 
\label{hodgsame}
\end{theorem}
\begin{proof}
The symmetry 
\[
g : [z_0 : z_1 : z_2 : z_3 : z_4]\mapsto  [\xi^{-1} z_0 : z_1 : z_2 : z_3 : z_4]. 
\]
descends to a symmetry of $\P^4/H$ and maps $\bar Z_t$ isomorphically to 
$\bar Z_{\xi t}$. By an argument analogous to the proof of 
Proposition~\ref{extends}, 
this isomorphism extends to 
an isomorphism $Z_t \rightarrow Z_{\xi t}$. This gives a diagram of 
polarized Hodge structures
\[\begin{array}{ccc}
H^3(Y_t, \Z)_{\rm free} & \stackrel{\pi^*}{\longrightarrow} & H^3(Z_t, \Z)_{\rm free} \\
&& {\Big\downarrow\rlap{$\vcenter{\hbox{$\scriptstyle\cong$}}$}}
 \\
H^3(Y_{\xi t}, \Z)_{\rm free} & \stackrel{\pi^*}{\longrightarrow} & H^3(Z_{\xi t}, \Z)_{\rm free}. 
\end{array}
\]
Comparing this with the action of $\Gamma$ on ${\mathcal D}_Y$ 
defined above gives the first statement. The second statement is immediate. 
\end{proof}

\begin{remark}\rm The proof of Lemma~\ref{modules} implies that
the spectral sequence
\[ E_2^{p,q}=H^p(K; H^q(Z, A)) \Rightarrow H^{p+q}(Y, A)
\]
degenerates at $E_2$ whenever 5 is invertible 
in $A$. In particular, there is an isomorphism of
polarized Hodge structures
\[ H^3 (Y_t, \Z[1/5]) \cong  H^3 (Y_{\xi t}, \Z[1/5]).\]
The problem is that $\aut(V_Y\otimes\Z[1/5], Q_Y)$ does not act properly 
discontinuously on ${\mathcal D}_Y$, so such a statement is weaker than the 
one proved above. On the other hand, it seems difficult to determine the 
precise behavior of the spectral sequence with $\Z$ coefficients, i.e. to 
compute the torsion in the cohomology of~$Y$. 
\end{remark}

\begin{remark}\rm The isomorphism of $\Q\,$-Hodge structures is due to 
Aspinwall and Morrison. They give a different proof coming from mirror 
symmetry which goes as follows. The mirror family ${\mathcal X}$ of $\by$ 
is the quotient of a suitable family of quintic hypersurfaces by the group 
$\langle\, g_1, g_3\,\rangle$. In particular, the antichiral ring of the 
central fibre $X_0$ of ${\mathcal X}$ with a choice of (complexified) 
K\"ahler class is isomorphic to the chiral ring of $Y_t$. On the other hand, 
the antichiral ring of $X_0$ can be shown to depend, via the mirror map, 
on $t^5$ only and not on $t$. Thus the varieties $Y_{\xi^i t}$ for 
$i=0,\ldots, 4$ have the same chiral ring, i.e. isomorphic rational Hodge 
structure. 
\end{remark}

\begin{remark}\rm Suppose that $Y_0$ is an $n$-fold, $G$ (a nontrivial quotient
of) the fundamental group $\pi_1(Y_0)$. Then there is an \'etale cover
$Z\rightarrow Y$; in fact there is a cover $Z_t\rightarrow Y_t$
for every deformation $Y_t$ of $Y_0$. The (primitive) cohomology 
$H^n_0(Z_t)$ becomes a $G$-representation, and in some cases one
can recover information about $Y_t$ from the pair 
\[(H^n_0(Z_t), \mbox{ action of } G).\]
A particular example of this construction 
is the theorem of Horikawa~\cite{horikawa}, giving a Torelli-type result for
Enriques surfaces using global Torelli for K3s. 
However, by Proposition~\ref{fundZ}, the threefold $Z_t$ under 
investigation is simply connected. On the other hand, as the proofs above
show, the Hodge structure on the middle-dimensional rational cohomology of the 
universal cover $Z_t$ contains no extra information, and it carries the 
trivial action of the fundamental group $\pi_1(Y_t)$. 
\label{nonact}
\end{remark}

\section{The automorphism group of the central fibre}
\label{am_central}

\begin{theorem} 
\label{autY} The automorphism groups of the varieties $Y_0$, $\bar Y_0$ are
\[\aut(Y_0)\cong\aut(\bar Y_0)\cong \langle G, g_4, g_5 \rangle / G,\] 
where 
\[
\begin{array}{rccl}
g_4: & [z_0 : z_1 : z_2 : z_3 : z_4] & \mapsto &  [z_0 :  z_1 :  z_2 : \xi^4 z_3 : \xi z_4]\\
g_5: & [z_0 : z_1 : z_2 : z_3 : z_4] & \mapsto &  [z_0 : z_2 : z_4 : z_1 : z_3].
\end{array}
\]
In particular, every automorphism of $\bar Y_0$ extends to an automorphism
on all (small) deformations $\bar Y_t$ of $\bar Y_0$. 
\label{explgroup}
\end{theorem}

\begin{proof} The proof of Theorem~\ref{autY} uses three Lemmas. 
The first one should certainly be well-known, 
but I could not find a suitable reference so I 
included a proof.

\begin{lemma}
Let 
\[
X = \left\{ \sum_{i=0}^n x_i^d =0 \right\} \subset \P^n_k
\]
be the Fermat hypersurface. 
Assume that $d\geq 3$, $n\geq 2$ and that $(n,d)\neq (2,3)$ or $(3,4)$.
Then
\[
\aut(X) \cong G_{n,d}, 
\]
where $G_{n,d}$ is the semi-direct product
$\Sigma_{n+1} \ltimes (\mu_d)^n$
of a symmetric group and a power of the group of $d$-th roots of unity. 
\label{am_f}
\end{lemma}
\begin{proof} For $n=2$, the result is proved in \cite{tzermias}. If $n\geq 3$ 
and $(n,d)\neq (3,4)$, then I first claim that every automorphism comes from a 
projective automorphism in the given embedding. If $n\geq 4$, Lefschetz 
implies $\pic(X)\cong\Z$ and then the claim is clear. If $n=3$ and $d\neq 4$ 
then the canonical class is (anti-)ample and this easily implies the claim 
again, see~\cite{matsumura_monsky}. 

Take an element $\sigma\in\aut(X)$ represented by an invertible matrix  
$A=(a_{ij})$. Apply $A$ to the equation of $X$ and consider the coefficients 
of $x_0^{d-1}x^{ }_1$, $x_0^{d-2}x_1^2$, and $x_0^{d-2}x^{ }_1 x^{ }_i$ for 
$i>1$. Their vanishing shows that the set of numbers 
\[\{a_{00}^{d-2}a_{01}^{ }, a_{10}^{d-2}a_{11}^{ },\ldots , a_{n0}^{d-2}a_{n1}^{ }\}\] 
solves the homogeneous system of equations given by the invertible matrix 
$A^T$. So all these quantities are zero. By symmetry, $a_{ij}a_{ik}=0$ 
whenever $j\ne k$. Hence $A$ has at most one non-zero entry in each row. 
Multiplying by a suitable element in $\Sigma_{n+1}$, $A$ can be brought into 
diagonal form, and then all its entries are $d$-th roots of unity.
\end{proof}

\begin{lemma} 
Let $\bar X$ be a canonical Calabi--Yau threefold with a finite number 
$m\geq 2$ of isolated $\frac{1}{5}(1,1,3)$ quotient singularities and Picard 
number one. Let $\pi:X\rightarrow \bar X$ be
the Calabi--Yau resolution. Then $\aut(X) \cong \aut(\bar X)$.
\label{cyisom}
\end{lemma}
\begin{proof} The Picard group of the resolution $X$ is
\[
\pic_\Q(X) \cong \Q H \oplus \Q E_1 \oplus \Q F_1 \oplus \ldots \oplus \Q E_m \oplus \Q F_m,
\]
where $H=\pi^*({\mathcal O}_{\bar X}(1))$ and $E_i$, $F_i$ are the classes
of the exceptional divisors as described in Proposition~\ref{unique}. The 
intersection numbers are as follows:

\begin{tabbing}
xxxx\=abcdefghxxxxxxxxxxxxxxxxxxxijklm\=\kill
\>$H^3 = d>0$\dotfill \>the degree of $\bar X$,  \\
\>$H\cdot E_i = H \cdot F_i=0 $\dotfill\>as $H$ is a pullback,\\
\>$E_i \cdot E_j = E_i \cdot F_j = F_i \cdot F_j = 0 $\dotfill\>unless $i=j$,\\
\>$E_i^3 = (K_{E_i})^2 =8$\dotfill\>as $E_i\cong {\mathbb F}_3$,\\
\>$F_i^3 = (K_{F_i})^2 =9$\dotfill\>as $F_i\cong \P^2$,\\
\>$E_i^2F_i = 1$,\\
\>$F_i^2E_i = -3.$
\end{tabbing}

\noindent Introducing the basis $H_0=H$, $H_{2i-1}=E_i+\frac{1}{3} F_i$, 
$H_{2i}=F_i$ of $\pic_\Q(X)$, the cubic form takes the shape
\[ \left(\sum_{i=0}^{2m} \alpha_i H_i\right)^3 = d \alpha_0^3 + 8 \frac{1}{3} \sum_{i=1}^m \alpha_{2i-1}^3 + 9 \sum_{i=1}^m \alpha_{2i}^3.  
\] 
Finally, the values of the second Chern class are 
\[  c_2(X)\cdot E_i = -4, \,
   c_2(X)\cdot F_i = -6, \,c_2(X)\cdot H = c\geq 0, 
\]
where the last inequality follows from a result of 
Miyaoka, \cite[Theorem 1.1]{miyaoka}. 

Let $\sigma\in\aut(X)$ be an automorphism. It acts via pullback on 
$\pic_\Q(X)$, fixing the cubic form together with the linear form given by cup
product with $c_2(X)$. I claim that the element $H_0=H$ of $\pic_\Q(X)$ must 
be fixed under the action. To see this, note that the cubic form has been 
manufactured to take the shape of the Fermat cubic. Every automorphism of
$\pic_\Q(X)$ must fix the associated (projectivized) hypersurface. The 
possible automorphisms are known from Lemma~\ref{am_f}. Moreover, in the 
present case, the multiplications by roots of unity are excluded since 
$\sigma$ must fix a {\it rational} vector space. The possible permutations are 
constrained by the fact that $c_2$ has to be fixed as well. As $c_2$ is 
negative on the $H_i$ for $i>0$ and non-negative on $H=H_0$, the latter is 
fixed and this proves the claim.

For large and divisible $m$, the divisor class $mH$ is base-point free and as 
the torsion in $\pic(X)$ is finite, is the unique representative of its 
numerical equivalence class. As $H\in\pic_\Q(X)$ is fixed by the induced 
action of $\sigma$, for large and divisible $m$ the space of sections of the
linear system $\mid mH \mid$ is also acted on by $\sigma$. In other words, 
the automorphism $\sigma$ descends to the image of the associated morphism 
which is exactly $\bar X$. 

For the converse, note that the quotient singularity $\frac{1}{5}(1,1,3)$ 
has a unique crepant resolution. Hence every automorphism 
$\bar\sigma\in\aut(\bar X)$ extends to a biregular automorphism 
$\sigma\in\aut(X)$ of the resolution. The Lemma follows.
\end{proof}

\begin{lemma} Let $X$ be a smooth algebraic variety with finite 
fundamental group~$F$. Let $Y$ be the universal cover of~$X$, a 
smooth algebraic variety with an action of~$F$ by automorphisms. Then
\[ \aut(X)\cong N_{\aut(Y)}(F)/F.
\]
\label{norm}
\end{lemma}
\begin{proof} Obvious. \end{proof}

To finish the proof of Theorem~\ref{autY}, let $Q^0_0$ be the open set of the 
Fermat quintic $Q_0$ on which the action of $G$ is free. Let $Y^0_0=Q_0^0/G$. 
There is a sequence of maps
\[ \aut(\bar Y_0) \hookrightarrow \aut(Y_0^0) \cong N_{\aut(Q_0^0)}(G)/G \cong N_{\aut(Q_0)}(G)/G. \]
The first isomorphism follows from Lemma~\ref{norm}. The second isomorphism 
uses $\aut(Q_0^0)\cong \aut(Q_0)$; here $\aut(Q_0^0)\subset \aut(Q_0)$ is 
proved by the argument used already in Proposition~\ref{fundZ} and the other 
direction is clear by Lemma~\ref{am_f}. 

On the other hand, by Lemma~\ref{am_f}, the automorphism group of $Q_0$ is 
the semi-direct product $G_{4,5}$ of the permutation and diagonal symmetries. 
Finding the normalizer of $G$ in $G_{4,5}$ is a finite search best done
using a computer; a short Mathematica routine computes this normalizer to be 
\[N_{\aut(Q_0)}(G)/G\cong\langle G, g_4, g_5 \rangle / G\]
with $g_4$, $g_5$ as in the statement of Theorem~\ref{autY}. So I obtain
\[\aut(\bar Y_0) \hookrightarrow \langle G, g_4, g_5 \rangle / G\]
and it is easy to see that this is in fact an isomorphism. Finally, by 
Lemma~\ref{cyisom}, $\aut(\bar Y_0)\cong\aut(Y_0)$. This proves the first 
statement. The second statement follows by inspection: every generator of the 
normalizer fixes $Q_t$. 
\end{proof}

\section{The proof of Theorem~\ref{am_main_theorem}}
\label{not_isomorphic!}

The proof is based on the following
rather standard result, a version of which was used in~\cite{tcex} already: 

\begin{theorem} Let $\bigx_i\rightarrow B$, $i=1,2$ be families of 
canonical Calabi--Yau varieties over a base scheme $B$, 
having simultaneous resolutions
$\by_i\rightarrow\bigx_i$ over $B$. 
Let $\bl_i$ be relatively ample
relative Cartier divisors on $\bigx_i$. Let ${\mathop{\rm Isom}}_B(\bigx_i, \bl_i)$ 
be the functor
\[
{\mathop{\rm Isom}}_B(\bigx_i, \bl_i): \underline{\rm Schemes} \rightarrow \underline{\rm Sets}\]
defined by
\[ {\mathop{\rm Isom}}_B(\bigx_i, \bl_i)(S)= \left\{\mbox{polarized }S\mbox{--isomorphisms } (\bigx_1)_S\rightarrow (\bigx_2)_S\right\},
\] 
where the pullback families $(\bigx_i)_S$ are polarized by
the relatively ample line bundles $(\bl_i)_S$. 
This functor is represented by a scheme ${\bf Isom}_B(\bigx_i, \bl_i)$,  
proper and unramified over $B$. 
\label{fp}
\end{theorem}
\begin{proof} By Grothen\-dieck's theory of the 
representability of Hilbert schemes and related functors, the above functor
is represented by a scheme ${\bf Isom}_B(\bigx_i, \bl_i)$, separated and
of finite type over $B$. The fact that the fibres have no infinitesimal 
automorphisms implies that ${\bf Isom}_B(\bigx_i, \bl_i)$ is unramified 
over $B$. Properness follows from the valuative 
criterion along the lines of~\cite[Proposition 4.4]{fantechi_pardini}; 
the existence of a simultaneous resolution is needed for this final step. 
\end{proof} 

\begin{theorem} Let $\by\rightarrow B$ be the family constructed in
Section~\ref{am_defs}, $\xi$ a primitive fifth root of unity. Then there is a 
Zariski dense subset $U\subset B$, such that the fibres $Y_t$ and $Y_{\xi t}$ 
are not isomorphic as algebraic varieties for $t\in U$.
\label{am_onestep_theorem}
\end{theorem}
\begin{proof} First I work with the singular family $\bar{\mathcal Y}$; 
for ease of notation, let $\bar{\mathcal Y}_1=\bar{\mathcal Y}$. 
Fixing an ample divisor $L$ on $\P_\Delta/K$ gives by restriction a 
relatively ample divisor $\bl$ on $\bar\by_1$. Let $\bl_1=\bl^{\otimes 5}$. 

Let $\gamma:B\rightarrow B$ be the map of the base which is multiplication by 
$\xi^{-1}$. Let $\bar\by_2\rightarrow B$ denote the pullback of 
$\bar\by_1\rightarrow B$ by $\gamma$. The family $\bar\by_2\rightarrow B$ 
is equipped with the relatively ample line bundle $\bl_2=\gamma^*(\bl_1)$ and 
its fibre over $t\in B$ is $\bar Y_{\xi t}$. 

\begin{lemma} Let $t\in B$, and let $\bar Y_{i,t}$ be the fibres of 
the two families polarized by the ample divisors $L_{i,t}$. Then every 
isomorphism \[\phi:\bar Y_{1,t}\stackrel{\sim}{\longrightarrow}\bar Y_{2,t}\] 
satisfies $\phi^*(L_{2,t})\sim L_{1,t}$.
\label{canonicalchoice}
\end{lemma}
\begin{proof} 
The fibres have Picard number one, and multiplication by five annihilates
every torsion element in their Picard groups. So the divisors 
$L_{i,t}$ are canonical elements of the respective Picard groups. 
The lemma follows.  
\end{proof}

Continuing the proof of Theorem~\ref{am_onestep_theorem}, consider the 
relative isomorphism scheme \[{\bf Isom}={\bf Isom}_B(\bar\by_i, \bl_i)\]  
together with the natural map ${\bf Isom}\rightarrow B$.
By Theorem~\ref{fp}, this map is proper, so its image
$V$ is a closed subvariety of the quasi-projective variety $B$. 

Assume first that $V=B$. Then 
${\bf Isom}$ has a component ${\bf I}$ with a 
surjective unramified map onto a Zariski neighbourhood of $0\in B$.
Now switch to the complex topology; let $\Delta$ be a disc in ${\bf I}$ 
mapping isomorphically onto a neighbourhood of $0\in B$.  
Consider the pullback families $\bar\by_{i,\Delta} \rightarrow \Delta$.  
By the definition of ${\bf I}$, these families are isomorphic under an 
isomorphism $\phi$ over $\Delta$. 

Consider the composition
\[\bar\by_{1,\Delta}\stackrel{\phi}{\longrightarrow} \bar\by_{2,\Delta}\stackrel{(\gamma^{-1})^*}{\longrightarrow}\bar\by_{1,\Delta}. \]
Its restriction to the central fibre $\bar Y_0$ is a polarized 
automorphism~$\sigma$.

By Proposition~\ref{betterprop}, $\bar\by_1\rightarrow \Delta$ 
is the universal deformation space of $\bar Y_0$ in the analytic category.  
The automorphism
$\sigma$~acts on the base of the deformation space by universality. 
This action equals the composite of the actions of $\phi$ and 
$(\gamma^{-1})^*$ on the base $\Delta$. However, $\phi$ is an isomorphism over 
$\Delta$, so the action of $\sigma$ on $\Delta$ is  multiplication 
by a primitive fifth root of unity, i.e. a rotation of the disc. 

On the other hand, by Theorem~\ref{autY}, the action of every automorphism
of $\bar Y_0$ on the base of the universal deformation space is {\it trivial}. 
Thus $\sigma$ cannot exist. So the assumption $V=B$ leads to a contradiction. 

Thus $V$ is a proper closed subset of $B$.  Let $U=B\setminus V$, 
a Zariski open subset of $B$. Over $t\in U$ 
the scheme ${\bf Isom}$ has no points. 
Using Lemma~\ref{canonicalchoice}, this implies that
for $t\in U$ there cannot exist any
isomorphism between $\bar Y_t$ and $\bar Y_{\xi t}$.

Finally, if $Y_t\cong Y_{\xi t}$ for
some $t\in B$, then an argument analogous to 
the proof of Lemma~\ref{cyisom} shows that the singular 
Calabi--Yau models $\bar Y_t$, $\bar Y_{\xi t}$ are also isomorphic. 
This concludes the proof of Theorem~\ref{am_onestep_theorem}.
\end{proof}

\noindent
Applying this theorem for $\xi^i$, $i=1,\ldots, 4$ and taking the
intersection of the resulting open sets concludes the
proof of Theorem~\ref{am_main_theorem} announced in the Introduction.

\begin{remark} \rm Theorem~\ref{am_main_theorem} is also argued for 
in the paper~\cite{aspinwall_morrison}. Aspinwall and Morrison write down a 
power series in the coordinate $t$ of the base $B$, following~\cite{bcov}, 
related to higher genus Gromov--Witten invariants of the family mirror family 
${\mathcal X}$. This series is a function of $t$ rather than $t^5$, and this 
is a strong indication of the validity of Theorem~\ref{am_main_theorem}. As a 
matter of fact, I believe that this is also an indication of the validity of 
Conjecture~\ref{notbir}. However, a solid mathematical definition, let alone 
computation, of this power series has not been given to date. 
\end{remark}

\bibliographystyle{amsplain}

\end{document}